\documentstyle[12pt]{article}
\normalbaselines
\begin{document}
\vbox {\vspace{6mm}}
\newcommand{\n}{{\bf n}}
\newcommand{\C}{{\bf C}}
\newcommand{\Z}{{\bf Z}}
\newcommand{\A}{{\bf A}}
\newcommand{\eps}{\varepsilon} 
\newcommand{\al}{\alpha}
\begin{center}
{\bf Fields of fractions of quantum solvable algebras}
\end{center}
\begin{center}
{\bf Alexander Panov}
\end{center}
\begin{center}
Department of Mathematical Sciences, Samara State Universety,
Russia, 443011, Samara, ul.Akad.Pavlova, 1\\
E-male: panov@info.ssu.samara.ru
\end{center}

We introduce the notion of pure $Q$-solvable algebra. 
The quantum matrices, quantum Weyl algebra, $U_q(\n)$ 
are the examples. It is proved that the skew field of fractions
of a pure $Q$-algebra $R$ is 
isomorphic to the skew field of twisted rational functions.
This is a quantum version of Gelfand-Kirillov conjecture for
solvable algebraic Lie algebras.

\begin{center}
{\bf 1.INTRODUCTION}
\end{center}

Consider a field $k$ of characteristic zero.
Let $Q=(q_{ij})$ be a $n\times n$-matrix with the entries $q_{ij}\in k*$
or $q_{ij}$ is an indeterminate, satisfying $q_{ij}q_{ji}= q_{ii}= 1$.
The algebra $k[Q]$ is generated by k and $q_{ij}^\pm$, $i,j=1,\ldots,n$.
Denote $k(Q)= Fract(k[Q])$ and $\Gamma $ is the subgroup in $k(Q)$
generated by $k^*$ and the entries of $Q$.
Throughout the paper $R$ is a $k$-algebra with a unit.\\
{\bf Definition 1.1.} We say that $R$ is a $Q$-algebra if $k[Q]$ is contained
in the center of $R$ (i.e. $R$ is a $k[Q]$-module).\\ 
{\bf Definition 1.2.} An algebra $R$ is $Q$-solvable if $R$ is a $Q$-algebra,
$R$ is freely generated as $k[Q]$-module by the elements
$x^{\overline{m}}=x_1^{m_1}x_2^{m_2}\cdots x_n^{m_n}$,
$\overline{m}=(m_1,m_2,\ldots,m_n)\in \Z_+^n$ and 
$$x_ix_j - q_{ij}x_jx_i = r_{ij}\eqno(1.1)$$ 
where $i<j$ and $r_{ij}$ is some element of the subalgebra
$R_{i+1}$ generated by $x_{i+1},\ldots,x_n$.\\
We remark that 
a $Q$-solvable algebra is an iterated skew polynomial extension of $k[Q]$,
considered in [GL].\\
A $Q$-solvable algebra admites the decreasing chain of subalgebras
$$R=R_1\supset R_2 \ldots \supset R_n\supset R_{n+1}=k[Q]$$ 
We consider the lexicographical order $<$ on $\Z_+^n$ such that
$e_1> e_2> \cdots> e_n$ where 
$e_i=(\delta_{i_1},\delta_{i,2},\ldots,\delta_{i,n})$. 
We denote by
$R(\overline{m})$ the $k[Q]$-submodule generated by the elements 
$\{ x^{\overline{p}} : \overline{p}\le\overline{m}\}.$
We define the degree for $a\in R$ as follows
$deg(a)= min \{\overline{m} : a\in R(\overline{m})\}$.

A $Q$-solvable algebra is a $\Z_+^n$-filtered algebra. Denote $A= gr(R)$.
The algebra $A$ is generated by $X_1,X_2,\ldots,X_n$ with
$X_iX_j= q_{ij} X_jX_i$. Thus $A$ is an algebra 
of twisted polynomials.
The $A$ is a Noetherian domain. Hence, every solvable $Q$-algebra $R$ is a 
Noetherian domain.

Let $R$ be an arbitrary $Q$-algebra. For every prime ideal $E\in Spec(k[Q])$
we denote $I_E= RE$, $R_E= R/I_E$ and $Z_E$ is the center of $R_E$. 
Consider the subset
$\Omega$ in $Spec(k[Q])$, consisting of the prime ideals $E$ such that
$R_E$ is a domain (this is always true if $R$ is $Q$-solvable) and 
$R_E$ is finitely generated as $Z_E$ module.\\
{\bf Definition 1.3.} We say that $R$ is a pure $Q$-algebra if $\Omega$
is dense in $Spec(k[Q])$ in Jacobson topology.\\
To check this condition it is sufficient to show that
$\Omega$ contains some family of ideals with zero intersection.
It is well known that the specialization of quantum algebras $q\mapsto\eps$ 
at roots
of unity leads to the algebras $R_\eps$ which is finitely generated over 
center.
This proves that all algebras of
the Examples are pure $Q$-algebras. 
Roughly speaking a $Q$-algebra is pure if it has sufficiently many
"good" specializations.

The Main Theorem of the paper asserts that the skew field of fractions 
of every pure $Q$-solvable
algebra is 
isomorphic to the skew field of twisted rational functions.
 
\begin{center}
{\bf 2. ON PURE QUANTUM ALGEBRAS}
\end{center}
We need the following\\
{\bf Definition 2.1}[AD1]. Let $R$ be an arbitrary Noetherian domain.
We say that $R$ is pure quantum if there is no embedding 
of the Weyl algebra $\A_1$ in the skew field $Fract(R).$\\
The Weyl algebra over field $k$ is generated by the elements $x$, $y$, related
by $xy-yx=1$. Recall that $\A_1$ has no nontrivial ideals. \\
{\bf Lemma 2.2.} 
Let $R$ be a Noetherian domain.
There exist the elements $x,y\in D=Fract(R)$ with $xy-yx=1$ if 
and only if there exist the non zero elements
$a,b,c,d,u,v,s,t,z,p,w,q$ in $R$, satisfying the conditions
$av=sd$, $ub=ct$, $zu=ps$, $vq=tw$,
$pabw-zcdq= zuvq$.\\
{\bf Proof.} See Lemma 2.9 in the paper [AD1].\\ 
{\bf Theorem 2.3.} Let $R$ be a Noetherian domain. Suppose that
$R$ is a pure $Q$-algebra and $k[Q]$ freely acts on $R$.
Then the algebra $R$ is pure quantum.\\
{\bf Proof}. Suppose that the Weyl algebra $\A_1$ is embedded in $D=Fract(R)$.
There exist the elements $x,y\in D$, related $xy-yx=1$. According to 
Lemma 2.2, there exist the elements $a,b,\ldots,q$, obeying the above
conditions. One of the following two cases may happen :\\
1) there exists the ideal $E\in\Omega$ such that all elements 
$a,b,\ldots,q$ don't belong to the ideal $I_E$ or\\
2) for every $E\in\Omega$ there exists an element in the set
$\{a,b,\ldots,q\}$ belonging to $I_E$.

Suppose that the case 2) holds. Then the product $\pi$ of the elements 
$\pi= ab\cdots q$ belongs to the ideal
$I_\Omega=\bigcap_{E\in\Omega}I_E$. We shall show that $I_\Omega=0$.
 
The $R$ contains a basis $\{v_\al\}$ over $k[Q]$. The element $\pi$ uniquely 
decomposes in the sum $\pi= \sum_\al v_\al c_\al$ with $c_\al\in C$.
The element $\pi$ lies in $I_E$ whenever $c_\al\in E$.
It implies that, if $\pi\in I_\Omega$ whenever 
$c_\al\in \bigcap_{E\in\Omega}E = 0$. We get $\pi=ab\cdots q=0$.
This contradicts to the statement that the elements 
$a,b,\ldots,q$ are non zero. 

Therefore, the case 1) takes place. Denote by $\overline a, \overline b,
\cdots, \overline q$ the images of $a,b,\ldots,q$ in $R_E$. The $R_E$ is
a Noetherian domain for every $E\in\Omega$. 
Denote $D_E= Fract(R_E)$. By Lemma 2.2, there exist the elements 
$\bar x, \bar y$, related 
$\bar x\bar y - \bar y\bar x = 1$, in 
$D_E$. Denote by $F_E$ the field $Fract (Z_E)$.
By definition of $\Omega$, $D_E$ has finite dimension over $F_E$.
Consider the mapping $x,y \mapsto \bar x,\bar y$ of the Weyl algebra
$\A_1(F_E)$ over the field $F_E$ to $D_E$. The Weyl algebra has no nontrivial
ideals; the kernel of the above mapping is trivial. We get the embedding of
the Weyl algebra in $D_E$. However, the $dim_{F_E} D_E< \infty$ and
  $dim_{F_E} \A_1(F_E)= \infty$. A contrudiction. The  $R$ is pure quantum.
$\Box$\\
{\bf Corollary 2.4.} If $R$ is pure $Q$-solvable, then $R$ 
is pure quantum.\\
{\bf Proof}. Every $Q$-solvable algebra is a Noetherian domain  and $k[Q]$
freely acts on $R$. $\Box$\\
Let $q$ be indeterminate and
let $q_{ij}= q^{s_{ij}}$ with $s_{ij}\in\Z$. In this case $k[Q]=k[q,q^{-1}]$.\\
{\bf Corollary 2.5.} Let $R$ be a $Q$-solvable algebra over $k[q, q^{-1}]$.
If there exists infinitely many specializations $\eps:k[q,q^{-1}]\mapsto k^*$ 
such that
$R_\eps$ is finitely generated as $Z_\eps$-module, 
then $R$ is pure quantum.\\
{\bf Proof} Is clear.

\begin{center}
{\bf 3.  ELEMENTS OF FINITE ADJOINT ACTION}
\end{center}
{\bf Definition 3.1.} We say that $x\in R$ is an
element of finite adjoint action  ( or $X$ is a FA-element) if 
$x$ is not a zero divisor and for every $y\in R$ there exists a polynomial
$f(t)= a_0t^N + a_1t^{N-1} +\cdots + a_N$, $a_0\ne 0$, $a_N\ne 0$ 
over $k$ such that 
$$ 
a_0x^Ny + a_1x^{N-1}yx +\cdots + a_Nyx^N = 0\eqno (3.1) 
$$ 
{\bf Proposition 3.2.} A FA-element generates an Ore set. \\
{\bf Proof.} We are going to prove that $\{x^i, i\in\Z_+\}$ is an Ore set. 
We rewrite (3.1) as 
$x(a_0x^{N-1}y +\cdots + a_{N-1}yx^{N-1}) = -a_Nyx^N.$
Denoting $b = -\frac{1}{a_N}(a_0x^{N-1}y +\cdots + a_{N-1}yx^{N-1})$,
we get $xb= yx^N$. The existence of $b_1$, 
obeying $b_1x= x^Ny$, is proved similarly.$\Box$

Denote by $R_x$ the localization of $R$
on the set generated by $x$. Denote $Ad_x(y)= xyx^{-1}$. 
One can rewrite (3.1) in the ring $R_x$ as 
$$ f(Ad_x)y= 0\eqno (3.2)$$
{\bf Proposition 3.3.} Let $R$ be a Noetherian domain and an algebra over
a field $F$. Suppose that 
$R$ is pure quantum  over $F$
and let $x$ be a FA-element in
$R$. We assert that, if $W$ is $Ad_x$-invariant  finite dimensional subspace
and the minimal polynolial of $Ad_x\vert_W$ decomposes into linear
factors over $F$, then $Ad_x$ is diagonable in $W$.\\
{\bf Proof}. Suppose that the Jordan form of $Ad_x$ is different from
diagonal. Then there exist the elements $u,v\in R_x$ such that
$Ad_xu=\al u$, $Ad_xv=u+\al v$, $\al\in F$.
Denote $Y=\al u^{-1}vx^{-1}$. 
The simple calculations 
$xYx^{-1}= Ad_xY= Ad_x(\al u^{-1}vx^{-1})=
\al (Ad_xu)^{-1}(Ad_xv)x^{-1}= \al\al^{-1}u^{-1}(u+\al v)x^{-1}= x^{-1} + Y$
yield $xY-Yx=1$. The Weyl algebra $\A_1$ has no nontrivial ideals;
the homomorphism $x,y\mapsto x,Y$ of the Weyl algenbra $\A_1$ to $D$
has no kernel.
The algebra, generated by $x$ and $Y$, is isomorphic to Weyl algebra.
This contradicts to the assumption that $R$ is pure quantum.
We conclude that the Jordan form of $Ad_x$ coinsides
with diagonal form. $\Box$  
 
\begin{center}
{\bf 4. MAIN THEOREM}
\end{center}
{\bf Definition 4.1.} An algebra $A_Q$ of twisted polynomials is generated
by the elements $y_1,\ldots,y_n$, satisfying the relations
$y_iy_j=q_{ij}y_jy_i$.

It is well known that $A_Q$ is a Noetherian domain. 
The skew field $Fract(A_Q)$ 
is called a skew field of twisted rational functions. \\
{\bf Definition 4.2.} We shall say that an element $x$ semicommutates 
with an element $y$ if $xy= q yx$ for some $q\in k[Q]$.\\
{\bf Lemma 4.3.} Let $R'$ be a $k(Q)$-algebra. 
We suppose that 
$R'$ is freely generated over $k(Q)$ by the elements
$x^{\overline{m}}=x_1^{m_1}\cdots x_i^{m_i}
x_{i+1}^{\pm m_{i+1}}\cdots x_n^{\pm m_n}$,
$\overline{m}=(m_1,m_2,\ldots,m_n)\in \Z_+^n$ and 
the generators obey  the defining relations:\\ 
i) $x_{i+1},\ldots,x_n$ semicommutates with all generators  $x_1,\ldots,x_n$;\\
ii) The generators $x_1,\ldots,x_n$ satisfy (1.1), i.e 
$$x_jx_s- q_{js}x_sx_j=r_{js}$$
where $j<s\le i$ and $r_{js}$ is some element in the subalgera $R'_{j+1}$
generated by 
$$x_{j+1}, \ldots , x_{i+1}^{\pm} , \ldots, x_n^\pm.$$
Then the element 
$x_i$ is a FA-element in $R'$ and ,futhermore, for every $y\in R'$
there exists a splitted polynomial 
$f(t)= (t-\beta_1)\cdots(t-\beta_N)$ with roots in $\Gamma$,
 satisfying (3.1).\\
{\bf Proof}. The algebra $R'$ admites the decreasing chain of subalgebras
$$R=R_1 \supset R_2\supset\ldots\supset R_{i+1}'.$$ 
The algebra $R'_{i+1}$ is an algebra of twisted Laurent polynomials.
The algebra $R$ is generated
as $R'_{i+1}$-module by the elements
$x^{\overline m}=x_1^{m_1}x_2^{m_2}\cdots x_i^{m_i}$
where $\overline{m}=(m_1,m_2,\ldots,m_n)\in \Z_+^i$.
We consider the lexicographical order $<$ on $\Z_+^i$ as in Intoduction.
We define the degree for $a\in R$ as follows
$deg(a)= min \{\overline{m}\in \Z_+^i: a\in R(\overline{m})\}$.
Denote $x=x_i$.\\
1) We shall prove at this step  that $x$ is an Ore element in $R'$.
Let $y\in R'$. We are going to prove that there exists $N\in\Z_+$
and $b_1$ such that $ x^Ny= b_1x$.
The existence of $M\in\Z_+$ and $b_2$ satisfying $yx^M= xb_2$
can de proved similarly.

We shall prove the statement by induction on $deg(y)$.
If $deg(y)=\bar{0}$, then $y\in R_{i+1}'.$ and it is easy.
We assume that the statement is proved for $y$ with $deg(y)< \overline{m}$.
Let $deg(y)=\overline{m}\ne 0.$ 
The conditions (1.1) implies that there exists $\beta\in\Gamma$ and $b\in R'$
such that 
$$xy =\beta yx +b,\eqno(4.1)$$
where $b=0$ or $ deg(b)<deg(y)$.
By the assumption of induction, there exist $l\in\Z_+$ and $b'\in R'$ such 
that $x^lb=b'x$. Therefore $x^{l+1}y= x^l(\beta yx +b) = 
x^l\beta yx + b'x = (x^l\beta y + b')x$. Denoting 
$N=l+1$ and
$ b_1=x^l\beta y + b'$,
we get $ x^Ny= b_1x$. \\
2) At the second step we shall prove that $x$ is a FA-element.
Consider the localization of $R'$ on $x$. Let $y\in R'$. 
The existence of polynomial $f(t)$, satisfying $f(Ad_x)y= 0$, will be proved 
by induction on $deg(y)$. 

If $deg(y)= 0$, then $y\in R_{i+1}'$ and the proof in easy.
Let $deg(y)=\overline{m}$ and suppose that the statement is proved for the 
smaller degrees. We use (4.1) and get :
$$(Ad_x- \beta)y= xyx^{-1}-\beta y = bx^{-1}\eqno(4.2)$$
where $b=0$ or $deg(b)< deg(y)$. By assumption, there exists a polynomial
$g(t)$ such that $ g(Ad_x)b=0$. Denote $f(t)=g(t)(t-\beta)$.
Then $f(Ad_x)y=g(Ad_x)(Ad_x-\beta)y=
g(Ad_x)(bx^{-1})= (g(Ad_x)b)x^{-1}= 0$.$\Box$\\
{\bf Main Theorem.} Suppose that $k$-algebra $R$ is pure $Q$-solvable. 
Then $Fract(R)$ is isomorphic $Fract (A)$ where $A=gr(R)$.
In particular, the skew field of fractions of $R$ is isomorphic 
to the skew field of twisted rational functions. \\
{\bf Proof}. We shall replace step by step the algebra $R$ 
by the algebra $R'$ with
generators, obeying the conditions i),ii) of Lemma 4.3 , such that 
$Fract(R)= Fract(R')$. Decreasing $i$, we get finally the
algebra $R'$ of twisted Laurent polynomials. This will conclude
the proof.

Denote by $S$ the Ore set 
generated by $\gamma-\gamma'$ where 
$\gamma,\gamma'\in \Gamma$, $\gamma\ne\gamma'$.
We consider the localisation 
$R_S $ of $R$ over $S$. This is our first step.

Suppose that we already have constructed the algebra $R'$ over $k[Q]_S$ 
with the generators  
$x_1,\ldots,x_i,x_{i+1}^\pm,\ldots,x_n^\pm$, satisfying i),ii), such that 
$Fract(R)= Fract(R')$. The algebra $R'$ is $\Z_+^{i}$-graded and 
$gr(R') = A'$ is generated by
$X_1,\ldots,X_i,X_{i+1}^\pm,\ldots,X_n^\pm$, $X_iX_j=q_{ij}X_jX_i$
where $q_{ij}$ are the same as in $(1.1)$.
As in Lemma 4.3 , we denote $x=x_i$.

According to Lemma 4.3, $x$ is a FA-element in 
$R'$. Denote by $R''$ the localization $R_x'$. 
We have $Fract(R')=Fract(R'')$.
We shall choose the elements $ x_1^{(0)},\ldots, x_{i-1}^{(0)}$ 
of $R''$ such that 
the $x_1^{(0)},\ldots, x_{i-1}^{(0)}, x_i^\pm, x_{i+1}^\pm,\ldots,x_n^\pm $
satisfy the conditions of Lemma 4.3.
As above $R_j'$ (resp.$R''_j$) is the subalgebra in $R'$,
generated by 
$x_{j},\ldots, x_{i-1},x_i, R_{i+1}' $
( resp. by $x_{j},\ldots, x_{i-1},x_i^\pm, R_{i+1}'$).
Let $y=x_j$ be one of generators $x_1,\ldots,x_{i-1}$ in $R_x'$.
By (1.1), 
$$xy - \gamma yx= r\in R_{j+1}'\eqno(4.3)$$ 
with $\gamma= q_{ij}\in \Gamma$.
Consider the minimal $Ad_x$-invariant $k(Q)$-subspace $W(y)$ in 
$R'\otimes k(Q)$.
Since $x$ is a FA-element, then $dim~W(y)<\infty.$
Let $\mu(t)$ be a minimal polynomial of $y$ with respect to $Ad_x$.
The $dim~W(y)$ equals to $deg(\mu)$.

The (4.3) implies $Ad_xy=\gamma y+rx^{-1}$. 
Then $0= \mu(Ad_x)y= \mu(\gamma)y + r_1$ where $r_1\in R''_{j+1}$.
Since $y=x_j\notin R_{j+1}'$, then $\mu(\gamma)=0$ and $\gamma$ is a root of 
$\mu(t)$. 

According to Lemma 4.3, $\mu(t)$ is decomposed into linear factors
with roots in $\Gamma$. 
By Proposition 3.3, $Ad_x$ is diagonal in $W(y)$.
Hence $\mu(t)=(t-\gamma_0)(t-\gamma_1)\cdots(t-\gamma_M)$,
with the different roots $\gamma=\gamma_0,\gamma_1,\ldots,\gamma_M\in \Gamma$.
Denote 
$$\mu_m(t)= (t-\gamma_0)\cdots\widehat{(t-\gamma_m)}\cdots(t-\gamma_M).$$
$$
x_j^{(m)}=\frac{1}    
{(\gamma_m-\gamma_0)\cdots(\gamma_m-\gamma_M)}\mu_m(Ad_x)y$$
Since $\gamma_m-\gamma_k\in S$, the element $x_j^{(m)}$ belongs to $R_x'$ 
One can decompose $y$ in the sum of eigenvectors 
$y = x_j= x_j^{(0)}+ x_j^{(1)}+\cdots+x_j^{(M)}$ 
where $Ad_xx_j^{(m)}= \gamma_mx_j^{(m)}$,
$m\in\overline{0,M}$. Here $M$ depends on $j$.
Observe that $x_j^{(m)}= const\cdot \mu_m(Ad_x)x_j$.
It follows that\\
1) the elements $x_i,x_{i+1},\ldots,x_n$ semicommutates with $x_j^{(m)}$
for all $j$ and $m$;\\
2)
if $m\ne 0$, then  $t-\gamma$ divides $\mu_m(t)$ and, therefore, 
$x_j^{(m)}\in R''_{j+1}$. We have $y=x_j=x_j^0~modR''_{j+1}$.

Consider $l$ such that $0\le j<l\le i-1$. 
As above one can decompose $x_l= x_l^{(0)}+ x_j^{(1)}+\cdots+x_j^{(P)}$ 
into a sum
of eigenvectors: $Ad_xx_l^{(p)}= \xi_px_l^{(p)}$, $p\in\overline{0,P}$.
The elements $x_j$ and $x_l$ satisfies (1.1):
 $$x_jx_l - q_{jl}x_lx_j= r_{j+1}\in R''_{j+1}\eqno(4.4)$$
We replace  $x_j$ and $x_l$ by their decompositions into sums
of eigenvectors in (4.4):
$$\begin{array}{c}(x_j^{(0)}+ x_j^{(1)}+\cdots + 
x_j^{(M)})(x_l^{(0)}+ x_l^{(1)}+\cdots+x_l^{(P)}) -\\ 
q_{jl}(x_l^{(0)}+ x_l^{(1)}+\cdots+x_l^{(P)})(x_j^{(0)}+ x_j^{(1)}
+\cdots+x_j^{(M)})
= r_{j+1}\in R''_{j+1}\end{array}\eqno(4.5)$$
Denote 
$$\Lambda= \{(m,p) : 0\le m\le M, 0\le p\le P, 
\gamma_m\xi_p=\gamma_0\xi_0\}$$ 
Remark that $\Lambda$ does not contain any pair $(0,p)$ besides $(0,0)$.

 The subralgebra $R''_{j+1}$ is $Ad_x$-invariant and $Ad_x$ is diagonable 
in every finite dimensional subspace.
Decompose $r_{j+1}$ into the sum of eigenvectors.
Let $r_{j+1}^0$ be the component with the eigenvalue $\xi_0\gamma_0$ 
in $r_{j+1}$. The element $r_{j+1}^0$ can be calculated via (4.5):
$$
x_j^{(0)}x_l^{(0)} - q_{jl}x_l^{(0)}x_j^{(0)} + \sum_{(m,p)\in\Lambda}
(x_j^{(m)}x_l^{(p)} - q_{jl}x_l^{(p)}x_j^{(m)}) = r_{j+1}^0$$
Since $x_j^{(m)}$ and $x_l^{(p)}$ lies in $R''_{j+1}$ for $m\ne 0$ and
all $p$, we have
$$
x_j^{(0)}x_l^{(0)} - q_{jl}x_l^{(0)}x_j^{(0)} = - \sum_{(m,p)\in\Lambda}
(x_j^{(m)}x_l^{(p)} - q_{jl}x_l^{(p)}x_j^{(m)}) + 
r_{j+1}^0\in R''_{j+1}\eqno (4.6)$$
Recall that $x_j= x_j^0modR''_{j+1}$. This proves that 
the elements
$x_1^0,\ldots, x_{i-1}^0, x_i^\pm, x_{i+1}^\pm,\ldots,x_n^\pm $
generate $R''=R_x'$. 
Moreover the elements 
$(x_1^0)^{m_1}\cdots (x_{i-1}^0)^{m_{i-1}}
x_i^{\pm m_{i+1}}\cdots x_n^{\pm m_n}$ with
 $(m_1,m_2,\ldots,m_n)\in \Z_+^n$ form a free basis of $R''$ over $k[Q]_S$.
The relations (4.6) is similar to (1.1).
It is clear that $Fract(R')= Fract(R'')$.

We replace $R'$ by $R''$ and begin the process from the beginning.
Finally we get the algebra $R^{(n)}$ of twisted Laurent polynomials
with $Fract(R^{(n)})= Fract(R)$. The $R^{(n)}$ is generated by 
some elemelents $y_1^\pm,\ldots,y_n^\pm$ with the relations
$y_iy_j=q_{ij}y_jy_i$. We obtain $Fract(R)= Fract(R^{(n)})= Fract(A)$.
 This proves that $Fract(R)$ is a skew field of twisted rational functions.
$\Box$\\
{\bf Definition 4.4}. As in Proof of Main Theorem $S$ is the Ore set 
generated by the  elements $\gamma-\gamma'$ for 
$\gamma,\gamma'\in \Gamma$ and $\gamma\ne\gamma'$.\\
Let $J$ be a completely prime ideal in $k[Q]$. 
We call $J$ an admissible ideal if it
does not intersect with $S$. 
Note that in the case $k[Q]=k[q,q^{-1}]$, $q$ is an indeterminate, 
the ideal $J=k[q,q^{-1}](q-\eps)$ 
is admissible if $\eps$ is not a root of unity.
The specialization $R\mapsto R_J$ (resp.$R\mapsto R_\eps$) 
will be called admissible specialization.\\
{\bf Corollary 4.5}.  Main Theorem is true for the admissible specializations
of pure $Q$-solvable algebras.\\
{\bf Proof.} The algebra $R$ is freely generated as $k[Q]$-module 
by the elements
$x^{\overline{m}}=x_1^{m_1}x_2^{m_2}\cdots x_n^{m_n}$, $m_i\in\Z_+$. 
In Main Theorem we construct step by step the knew set of generators which are
defined over $k[Q]_S$. Finally we get the algebra $R^{(n)}$ of Laurent 
twisted polynomials with the semicommutating generators
$y_1^\pm,\ldots ,y_n^\pm$. The $R^{(n)}$ is freely generated as 
$k[Q]_S$-module by the elements
$y^{\overline{m}}=y_1^{m_1}y_2^{m_2}\cdots y_n^{m_n}$, $m_i\in\Z$ 
The ideal $J$ does not intersect with $S$. 
The algebra $R^{(n)}_J$ is an algebra of Laurent twisted polynomials
and $Fract(R_J)=Fract(R^{(n)}_J)$.$\Box$\\
\begin{center}
{\bf 5. EXAMPLES}
\end{center}
{\bf Example 1. Quantum matrices}.\\
Let $P,Q\in Mat_n(k)$ be matrices such that
$p_{ij}q_{ij}= c^{sgn(j-i)}$, $p_{ij}p_{ji}=q_{ij}q_{ji}=p_{ii}=q_{ii}=1$.
The algebra $M_{P,Q,c}(n)$ of regular functions on quantum matrices
is generated by the elements
$\{a_{ti}\}_{t,i=1}^n$ with the relations 
$$a_{ti}a_{sj}- q_{ts}q_{ij}^{-1}a_{sj}a_{ti}= 
(p_{ts}^{-1}- q_{ij})a_{si}a_{tj}$$
for $i<j$, $t<s$ or $i>j$, $t>s$ and
$$ a_{ti}a_{sj}= p_{ts}^{-1}q_{ij}^{-1}a_{sj}a_{ti} $$
in all other cases (see [MP]). Denote $x_{(i-1)n+j}=a_{ij}$ and $R_s$ is the 
algebra generated by $x_t$, $t=1,\ldots,s$. Consider the filtraton:
$$R=R_{n^2}\supset \ldots \supset R_1\supset R_{0}=k[Q]$$
The $M_{P,Q,c}(n)$ is a $Q$-solvable with respect to this filtration.
Consider the specialzations $p_{ij},c\mapsto \eps_{ij}, \eps_o$ 
at roots of unity $\eps_{ij}^{l_{ij}}=1$, $\eps_0^{l_0}=1$.
The elements $a_{ij}^l$ with $l=l_0\prod l_{ij}$ lie in the center
of $R_\eps$. The algebra $R_\eps$ is finitely generated over the center.
Therefore $M_{P,Q,c}(n)$ is pure $Q$-algebra.
Well known that monomials of $a_{ij}$
form a $k$-basis (see [D], 2.4).
The Main Theorem asserts that the skew field  $Fract(M_{P,Q,c}(n))$ 
is isomorphic to the skew field of fractions of the algebra of twisted
polynomials generated by 
$\{y_{ti}\}_{t,i=1}^n$ with the relations 
$y_{ti}y_{sj}= q_{ts}q_{ij}^{-1}y_{sj}y_{ti}$
for $i<j$, $t<s$ or $i>j$, $t>s$ and
$ y_{ti}y_{sj}= p_{ts}^{-1}q_{ij}^{-1}y_{sj}y_{ti}$
in all other cases (see [P],[MP],[C]).\\
{\bf Example 2. Quantum Weyl algebra}.\\
 The matrices $P,Q$ and $c$ as in Quantum matrices.
A Quantum Weyl algebra $A_{P,Q,c}(n)$ 
is generated by $x_i$,$y_i$, $i=\overline{1,n}$ with relations\\ 
$x_ix_j=q_{ij}^{-1}x_jx_i$;\\ 
$y_iy_j=p_{ij}y_jy_i$ ;\\
$y_ix_j=p_{ji}x_jy_i$ for $i<j$;\\ 
$y_ix_j=q_{ij}x_jy_i$ for $i>j$;\\
$y_ix_i=1+c^{-1}x_iy_i+(c^{-1}-1)\sum_{\al=i+1}^nx_\al y_\al$ ; 
see [AD2] and ([D],4.2).
\\
The Quantum Weyl algebra is an iterated skew polynomial extension of $k[Q]$ 
([GL]).
Therefore the  monomials  of $y_1,\ldots,y_n,x_1,\ldots,x_n$ form $k[Q]$-basis.
The Quantum Weyl algebra is $Q$-solvable.
Consider the specialzations $p_{i},c\mapsto \eps_{i}, \eps_o$ 
into roots of unities $\eps_{i}^{l_{i}}=1$, $\eps_0^{l_0}=1$.
The elements $x_{i}^l$ and $y_i^l$ with $l=l_0\prod l_{i}$ lie in the center
of $A_\eps= (A_{P,Q,c}(n))_\eps$. 
The algebra $A_\eps$ is finitely generated over the center.
Therefore $A_{P,Q,c}(n)$ is pure $Q$-algebra. The Main Theorem implies
that the skew field  $Fract(A_{P,Q,c}(n))$ 
is isomorphic to the skew field of fractions for the twisted
polynomial algebra, generated by $X_i$,$Y_i$, $i=\overline{1,n}$ 
with relations\\ 
$X_iX_j=q_{ij}^{-1}X_jX_i$;\\ 
$Y_iY_j=p_{ij}Y_jY_i$ ;\\
$Y_iX_j=p_{ji}X_jY_i$ for $i<j$;\\ 
$Y_iX_j=q_{ij}X_jY_i$ for $i>j$;\\
$Y_iX_i=c^{-1}X_iY_i$ (see [AD2]).\\
{\bf Example 3. $U_q(\n^+)$}\\
The quantum universal enveloping algebra $R= U_q(\n^+)$ 
for the upper nilpotent subalgebra $\n^+$
is generated over $k[q,q^{-1}, (q^{d_i}-q^{-d_i})^{-1}]$
by $E_i$, $i=\overline{1,n}$ with the quantum 
Chevallley-Serre relations.
Fix a reduced expression $w_0=s_{i_1}\ldots s_{i_N}$ of the longest element
in the Weyl group W. Consider the following convex ordering 
$$\beta_1=\al_{i_1}, \beta_2=s_{i_1}(\al_2),\ldots, 
\beta_N= s_{i_1}\ldots s_{i_{N-1}}(\al_N)$$
in the set $ \Delta^+$ of positive roots.
Introduce the corresponding $root~vectors$ for $s=\overline{1,N}$ ([L]):
  $$ E_{\beta_s}=T_{i_1}\cdots T_{i_{s-1}}E_{i_s}.$$
For $m=(m_1,\ldots,m_N)$  let 
$E^m= E_{\beta_1}^{m_1}\cdots E_{\beta_N}^{m_N}.$
The elements $E^m$, $m\in \Z_+^N$ form a basis of $U_q(\n^+)$
over $k[q,q^{-1}, (q^{d_i}-q^{-d_i})^{-1}]$ (see [L]). 
The relations between $E_{\beta_i}$, $i=\overline{1,N}$
as follows [LS]:
$$E_{\beta_i} E_{\beta_j}- q^{-(\beta_i,\beta_j)}
E_{\beta_j}E_{\beta_i}= \sum_{m\in\Z_+^N}c_m E^m,\eqno (5.1)$$
where $i<j$, $c_m\in k[q,q^{-1}]$ and $c_m\ne 0$ only 
when $m=(m_1,\ldots,m_N)$ is such that $m_s=0$ for $s\le i$ and $s\ge j$.
The element in the right side of (5.1) belongs to the subalgebra
$R_{i+1}$ generated by $E_{\beta_{i+1}},\ldots,E_{\beta_N}$
The relations (5.1) have the form (1.1). The algebra $U_q(\n^+)$ is
$Q$-solvable. The elements $E_\al^l$, $\al\in\Delta^+$ lie
in the center of $R_\eps$ if
$l$ is relatively prime to all $d_i$, $l> max\{d_i\}$ and $\eps$
is a primitive $l$-th root of unity [DCK]. Hence $U_q(\n^+)$ is pure
$Q$-algebra.
The Main Theorem asserts that the skew field  
$Fract(U_q(\n^+)) $
is isomorphic to the skew field of fractions of the algebra of twisted
polynomials generated by 
$e_{\beta_i}$, $i=\overline{1,N}$
with the relations $e_{\beta_i} e_{\beta_j}= q^{-(\beta_i,\beta_j)}
e_{\beta_j}e_{\beta_i}$ , $i<j$ (see [J]). By Corollary 4.5, the Main Theorem
remains true in the case $q\in k^*$ and $q$ is not a root of unity.\\

{\bf Acknowledgments}. I am very grateful to J.Alev and F.Dumas for
support. I would like to acknowledge P.Caldero for the following 
information: the Main Theorem of the article for $U_q(\n)$ 
and in some other cases
can be derived from Prop.2.2, Prop.2.3 of the paper [C1] and 
also from paragraph 3 of [C2].

\end{document}